\documentclass[10pt]{elsart}

\usepackage{amsmath}
\usepackage{amssymb}
\usepackage{hyperref}

\newcommand {\lG}{\log\Gamma}
\newcommand {\lGs}{\log^2\Gamma}
\newcommand {\lB}{\log G}
\newcommand {\Cl}{\mathrm{Cl}}

\usepackage[a4paper, left=1in,right=1in,top=1in,bottom=1in]{geometry}%====================================================================

\begin{document}

\begin{frontmatter}

 \title{The Fourier series of the log-Barnes function}
 \author[mezo]{Istv\'an Mez\H{o}}
 \ead{istvanmezo81@gmail.com}
 \address{Department of Mathematics,\\Nanjing University of Information Science and Technology,\\No.219 Ningliu Rd, Pukou, Nanjing, Jiangsu, P. R. China}
 \thanks[mezo]{The research of Istv\'an Mez\H{o} was supported by the Scientific Research Foundation of Nanjing University of Information Science \& Technology, and The Startup Foundation for Introducing Talent of NUIST. Project no.: S8113062001}

\begin{abstract}
In this paper we determine the Fourier series expansion of the log-Barnes function. This is the analogue of the classical result of Kummer and Malmsten.

Applying this expansion we get some integrals similar to the Espinosa-Moll log-Gamma integrals with respect to $\log G$. During the course of the paper some interesting log-Gamma integrals and a log-harmonic sum are evaluated.
\end{abstract}

\end{frontmatter}

\section{Raabe's formula and its extensions}

In 1840 J. L. Raabe \cite{Raabe} proved that
\[\int_0^1\lG(x+t)dx=\log\sqrt{2\pi}+t\log t-t\quad(t\ge0)\]
with the limiting case
\begin{equation}
\int_0^1\lG(x)dx=\log\sqrt{2\pi}.\label{RaabeOrig}
\end{equation}
Here $\Gamma$ is the Euler Gamma function
\[\Gamma(x)=\int_0^{\infty}u^{x-1}e^{-u}du\quad(\Re(x)>0).\]
During the studying of some definite integrals of the Hurwitz zeta function O. Espinosa and V. Moll \cite{EM1,EM2} proved an extension of the Raabe formula \eqref{RaabeOrig}:
\[\int_0^1\lGs(x)dx=\log\sqrt{2\pi}=\frac{\gamma^2}{12}+\frac{\pi^2}{18}+\frac13\gamma L_1+\frac43L_1^2-(\gamma+2L_1)\frac{\zeta'(2)}{\pi^2}+\frac{\zeta''(2)}{2\pi^2},\]
where
\[L_1=\int_0^1\lG(x)dx=\log\sqrt{2\pi},\]
$\gamma$ is the Euler constant and $\zeta$ is the Riemann zeta function. With this discovery a number of researchers started to look for exact values of the integrals of the powers of the log-Gamma function $\lG$. Introducing the notation
\[L_n=\int_0^1\log^n\Gamma(x)dx,\]
Amdeberhan et al. \cite{ACEKMM} reproduced the result of Raabe by an elementary argument, and reproved the Espinosa-Moll expression for $L_2$ using Parseval's identity and Kummer's Fourier series\footnote{I. Blagouchine \cite{Blagouchine} discovered that this expansion was known 5 years before Kummer by Malmsten and his colleagues.}\cite[entries 6.443.1, 6.443.3]{GR} or \cite[p. 29]{SCbook}
\begin{equation}
\lG(x)=\frac12\log(2\pi)+\sum_{n=1}^\infty\frac{C+\log n}{\pi n}\sin(2n\pi x)+\sum_{n=1}^\infty\frac{1}{2n}\cos(2n\pi x)\label{Kummer}
\end{equation}
when $0<x<1$ and $C=\log(2\pi)+\gamma$.

In \cite{ACEKMM} it was also proved that
\[L_n\sim n!\quad\mbox{as}\quad n\to\infty,\]
and the authors called for expressions of $L_n$ in general.

One year later Bailey and his co-authors \cite{BBB} found an expression for $L_3$ and $L_4$ in terms of Tornheim-Witten zeta functions and its derivatives. They remarked that in principle their approach might be used to find $L_n$ in general.

In our paper we would like to study analogous questions with respect to the Barnes $G$ function, which is a higher order analogue of the Gamma function.

\section{The Barnes $G$ function}

The Barnes $G$ function is defined by the functional equation
\begin{equation}
G(z+1)=\Gamma(z)G(z)\quad(z\in\mathbb C\setminus\{0,-1,-2,\dots\})\label{GzG}
\end{equation}
with the normalization $G(1)=1$. The origin of the Barnes function is in the works of Barnes \cite{Barnes1,Barnes2,Barnes3}. Later he defined a sequence of functions (the multiple gammas) as higher order analogues of the $\Gamma$ and $G$ functions \cite{Barnes4}. However, before the papers of Barnes, Kinkelin \cite{Kinkelin} defined a function similar to $G$ and applied it to evaluate some log-trigonometric integrals.

The $G$ function has applications in the theory of infinite sums involving the zeta function \cite{CCS,CSA}, in the theory of the Polygamma \cite{Adamchik,EMpolygamma}, in the evaluation of Laplacians \cite{Vardi} and in the evaluation of several integrals containing Gamma function \cite{CSA}. The asymptotic behavior of $G$ was discussed by C.-P. Chen \cite{Chen}, the asymptotics of $\lB$ is contained in \cite{Voros}, and other computational aspects were studied by Adamchik \cite{Adamchikcomp}. Additional applications are mentioned in \cite{Adamchikarxiv}. A good resource with many important properties of $G$ and $\lB$ is the book of Srivastava  and Choi \cite{SCbook}. Citations to applications can also be found in this book on p. 38.

The purpose of this paper is threefold: we study the integral
\[G_2=\int_0^1\log^2G(x)dx,\]
find the Fourier expansion of $\lB$ and find the asymptotics of $G_n$ in general, where
\[G_n=\int_0^1\log^nG(x)dx\quad(n\ge1).\]

\section{The Fourier series of $x\lG(x)$ and $\lB(x)$}

D. Connon \cite{Connon} and before Koyama and Kurokawa \cite[Theorem 1]{KK} found a Fourier-like expansion of $\lB$ as
\begin{equation}
\lB(x)=\zeta'(-1)-\frac{1}{4\pi}\sum_{n=1}^\infty\frac{1}{n^2}\sin(2n\pi x)+\frac1{2\pi^2}\left(C-\frac32\right)\sum_{n=1}^\infty\frac{1}{n^2}\cos(2n\pi x)+\label{Connon}
\end{equation}
\[\frac{1}{2\pi^2}\sum_{n=1}^\infty\frac{\log n}{n^2}\cos(2n\pi x)+\frac14B_2(x)+(x-1)\lG(x),\]
where $B_2(x)=x^2-x+\frac{1}{6}$ is the second Bernoulli polynomial and, as before, $C=\log(2\pi)+\gamma$. (We must note that The Koyama-Kurokawa Fourier series was expressed for the second multiple Gamma function $\Gamma_2$.)

One can then deduce the pure Fourier representation of $\lB$ by considering this expansion for $x\lG(x)$. This is not entirely elementary, so we now deduce the Fourier series of $x\lG(x)$. To fix the notation we consider the Fourier coefficients of a function $f(x)$ as
\begin{align*}
a_n&=2\int_0^1f(x)\cos(2n\pi x)dx\quad(n\ge0),\\
b_n&=2\int_0^1f(x)\sin(2n\pi x)dx\quad(n\ge1).
\end{align*}
Then
\[f(x)=\frac12a_0+\sum_{n=1}^\infty a_n\cos(2n\pi x)+\sum_{n=1}^\infty b_n\sin(2n\pi x)\quad(0<x<1).\]
The generalization of Parseval's identity will be important to us. This states \cite[p. 343]{Bartle} that if $(a_n,b_n)$ and $(\alpha_n,\beta_n)$ are the Fourier coefficients of $f$ and $g$, respectively, then
\begin{equation}
2\int_0^1f(x)g(x)dx=\frac12a_0\alpha_0+\sum_{n=1}^\infty(a_n\alpha_n+b_n\beta_n).\label{Parsgen}
\end{equation}

One more ingredient we will need in the course of the proof, which is the following easily provable lemma.

\begin{lem}\label{Lem1}The series identities
\[\sum_{m=1,\;m\neq n}^\infty\frac{1}{m^2-n^2}=\frac34\frac{1}{n^2}\quad(n\ge1),\]
\[\sum_{m=1,\;m\neq n}^\infty\frac{1}{m(m^2-n^2)}=\frac{5}{4n^3}-\frac{H_n}{n^2}\quad(n\ge1)\]
hold true. Here
\[H_n=\frac11+\frac12+\cdots+\frac1n\]
is the $n$th harmonic number.

\end{lem}

Having these tools we are going to present the Fourier series of $x\lG(x)$.
\begin{thm}\label{xlogGammaFourier}For any $0<x<1$
\[x\lG(x)=\frac{a_0}{2}+\sum_{n=1}^\infty a_n\cos(2n\pi x)+\sum_{n=1}^\infty b_n\sin(2n\pi x),\]
where
\begin{align*}
a_0&=\frac12\log\left(\frac{2\pi}{A^4}\right),\\
a_n&=\frac{1}{4n}-\frac{C}{\pi^2n^2}-\frac{\log n}{4\pi^2n^2}-\frac{1}{\pi^2}T_n\quad(n\ge1),\\
b_n&=\frac{\gamma+\log(n)-H_n}{2\pi n}+\frac{1}{2\pi n^2}\quad(n\ge1)
\end{align*}
with
\begin{equation}
T_n=\sum_{i=2\;(i\neq n)}^\infty\frac{\log(i)}{i^2-n^2}.\label{Tndef}
\end{equation}
Here
\[A=e^{1/12-\zeta'(-1)}\]
is the Glaisher-Kinkelin constant.
\end{thm}

\textit{Proof.} By definition,
\[a_n=2\int_0^1x\lG(x)\cos(2n\pi x)dx\quad(n\ge0).\]
To determine $a_n$ we can apply \eqref{Parsgen} with $g_n(x)=x\cos(2n\pi x)$ and $f(x)=\lG(x)$. Since the Fourier coefficients of $\lG(x)$ are known by \eqref{Kummer}, we can concentrate on the $(\alpha_{n,m},\beta_{n,m})$ Fourier coefficients of $g_n(x)$. It is not hard to see that
\[\alpha_{0,0}=1,\quad\alpha_{n,n}=\frac12\;(n\ge1),\]
and the other $\alpha$'s are zero. Moreover,
\[\beta_{m,m}=-\frac{1}{4\pi m}\;(m>0),\quad\beta_{n,m}=\frac{m}{\pi(n^2-m^2)}\;(n\neq m).\]
Thus, by \eqref{Parsgen}
\[a_0=\int_0^1x\lG(x)dx=\frac{1}{2}\log(2\pi)+\sum_{m=1}^\infty\left(\frac{1}{2m}\cdot0+\frac{B+\log m}{\pi m}\frac{m}{\pi(0^2-m^2)}\right)=\]
\[\frac12\log(2\pi)-\frac{1}{\pi^2}\sum_{m=1}^\infty\frac{A+\log m}{m^2}.\]
The sum on $m$ equals to $A\pi^2/6-\zeta'(2)$. The derivative of zeta at $x=2$ can be expressed by $\zeta'(-1)$ using the functional equation of $\zeta$. Thus
\[a_0=\frac12\log\left(\frac{2\pi}{A^4}\right).\]
We continue with $a_n$:
\[a_n=\sum_{m=1,\;m\neq n}^\infty\left(\frac{1}{2m}\cdot0+\frac{C+\log m}{\pi m}\frac{m}{\pi(n^2-m^2)}\right)+\left(\frac{1}{2m}\cdot\frac12+\frac{C+\log n}{\pi n}\frac{-1}{4\pi n}\right)=\]
\[\frac{1}{4n}-\frac{C+\log n}{4\pi^2 n^2}-\frac{1}{\pi^2}\sum_{m=1,\;m\neq n}^\infty\frac{C+\log m}{m^2-n^2}.\]
Simplifying this expression by using Lemma \ref{Lem1}, and after defining $T_n$ as above, we arrive at the expression for $a_n$ presented in the statement.

The expression for $b_n$ can be given by considering the $(\alpha_{n,m},\beta_{n,m})$ Fourier expression of $x\sin(2n\pi x)$ together with the second statement of \ref{Lem1}. The details are left to the reader.

%These are as follows:
%\[\alpha_{n,0}=-\frac{1}{\pi n}\;(n\ge1),\quad\alpha_{n,n}=-\frac{1}{4\pi n}\;(n\ge1),\quad\alpha_{n,m}=-\frac{n}{\pi(n^2-m^2)}\;(n\neq m).\]
%All the $\beta$'s are zero but $\beta_{n,n}=\frac12$. Hence
%\[b_n=\frac12\log(2\pi)(-\frac{1}{\pi n})+\sum_{m=1,\;m\neq n}^\infty\left(\frac{1}{2m}\frac{-n}{\pi(n^2-m^2)}+\frac{C+\log m}{\pi m}0\right)+\]
%\[\underbrace{+\frac{1}{2n}(-\frac{1}{4\pi n})+\frac{C+\log n}{\pi n}\frac12}_{\textrm{when }m=n}=\]
%\[-\frac{\log(2\pi)}{2\pi n}-\frac{1}{8\pi n^2}+\frac{C+\log n}{2\pi n}+\frac{n}{2\pi}\sum_{m=1,\;m\neq n}^\infty\frac{1}{m(m^2-n^2)}.\]
%Applying the second identity in Lemma \ref{Lem1} and taking into account the definition of $C$, we are done.
\hfil\qed

The pure Fourier expansion of the log-Barnes function now follows.

\begin{thm}\label{FourierLogG}For any $0<x<1$
\[\lB(x)=a_0+\sum_{n=1}^\infty a_n\cos(2n\pi x)+\sum_{n=1}^\infty b_n\sin(2n\pi x),\]
where
\begin{align*}
a_0&=\frac{1}{12}-2\log(A)-\frac14\log(2\pi),\\
a_n&=\frac{1}{2\pi^2n^2}\left(\frac{\log n}{2}-C-1\right)-\frac{1}{4 n}-\frac{T_n}{\pi ^2}\quad(n\ge1),\\
b_n&=\frac{1}{2\pi n}\left(\frac{1}{2n}-\gamma-\log (4 \pi^2n )-H_n\right)\quad(n\ge1)
\end{align*}
with $T_n$ is as in \eqref{Tndef}.
\end{thm}

\textit{Proof.} Theorem \ref{xlogGammaFourier} and Connon's formula \eqref{Connon} directly gives the statement if we apply the expansion \cite[p. 338]{Apostol}
\[B_2(x)=\frac{1}{\pi^2}\sum_{n=1}^\infty\frac{\cos(2n\pi x)}{n^2}.\]
\hfill\qed

\section{The integral $G_2$}

It was known even by Barnes that
\[G_1=\int_0^1\lB(x)dx=\frac{1}{12}-\frac14\log(2\pi)-2\log A.\]
See \cite{SCbook} for the details on $G_1$. Our aim is to find
\[G_2=\int_0^1\log^2G(x)dx.\]

In principle, to find $G_2$ we could use Parseval's identity together with Theorem \ref{FourierLogG} but there are so many terms under the square that this work would be very long. Instead, we will use the trivial identity
\begin{equation}
\log^2 G(x)=\frac14\left(\log^2G(x)G(1-x)+\log^2\frac{G(x)}{G(1-x)}\right),\label{G2expr}
\end{equation}
and determine the integrals
\[\int_0^1\log^2G(x)G(1-x)dx,\quad\mbox{and}\quad\int_0^1\log^2\frac{G(x)}{G(1-x)}dx\]
which are interesting in themselves. In the evaluation of these integrals three quantities appear which we could not trace back to more elementary expressions. These are the following.
\begin{align*}
\zeta_H'(2)&=\sum_{n=2}^\infty\frac{H_n\log n}{n^2}\approx2.623865966,\\
U&=\int_0^1Z^2(t)dt\approx0.4785935,\\
V&=\sum_{n=1}^\infty\frac{\zeta(2n+1)}{n+1}E_{2n+2}\approx0.055645894.
\end{align*}
Here
\[Z(t)=\sum_{n=1}^\infty\frac{\zeta(2n+1)}{n+1}t^{2n+2}\quad(|t|<1),\]
while $E_n$ is defined in \eqref{Endef}. The notation $\zeta_H'(2)$ will be explained later in a separate subsection.

The Clausen function is an odd function with $1/n^2$ sine-coefficients.
During the course of the below proof we shall need the Fourier coefficients of $x\Cl_2(x)$, where
\[\Cl_2(x)=\sum_{n=1}^\infty\frac{\sin(nx)}{n^2}\]
is the Clausen function. The following lemma contains this information. (The proof is omitted.)

\begin{lem}\label{LemCl2}The Fourier-cosine coefficients of $x\Cl_2(2\pi x)$ are
\[a_0=-\frac{\zeta(3)}{\pi},\quad a_n=\frac{1}{\pi}\left(\frac{H_n}{n^2}-\frac{3}{2n^3}\right)\;(n>0).\]
\end{lem}

%\textit{Proof.} The zeroth coefficient $a_0$ equals to
%\[a_0=2\int_0^1x\Cl_2(2\pi x)dx=-\frac{\zeta(3)}{\pi},\]
%as one can evaluate by partial integration, for example. In general,
%\[a_n=2\int_0^1x\Cl_2(2\pi x)\cos(2n\pi x)dx=2\Im\int_0^1x\left(\sum_{k=1}^\infty\frac{e^{2\pi ikx}}{k^2}\right)\cos(2n\pi x)dx=\]
%\[2\sum_{k=1}^\infty\frac{1}{k^2}\Im\int_0^1xe^{2\pi ikx}\cos(2n\pi x)dx=\]
%\[2\sum_{k=1\;k\neq n}^\infty\frac{1}{k^2}\left(\frac{-k}{2k^2\pi-2n^2\pi}\right)+\underbrace{2\frac{1}{n^2}\cdot\frac{-1}{8n\pi}}_{\textrm{when }k=n}=\frac{1}{\pi}\sum_{k=1\;k\neq n}^\infty\frac{1}{k(n^2-k^2)}-\frac{1}{4n^3\pi}.\]
%Applying Lemma \ref{Lem1} we are done after a simplification.
%\hfill\qed

We now evaluate the following log-Barnes integrals.

\begin{thm}We have that
\begin{equation}
\int_0^1\log^2\frac{G(x)}{G(1-x)}dx=2\log C(\gamma+2\log(2\pi))+\frac{\zeta ''(2)}{2 \pi ^2}+\frac{\zeta '(3)}{2 \pi ^2}+\label{GperG}
\end{equation}
\[\frac{1}{\pi ^2}\zeta_H'(2)+\frac{3\zeta (3)}{\pi ^2}\left(\frac{\gamma}{2}+\log (2\pi)\right)+\frac{13 \pi ^2}{720}-\frac{\gamma ^2}{12}-\frac{\gamma}{6}\left(\log2-\log (\pi )\right),\]
and
\[\int_0^1\log^2G(x)G(1-x)dx=\]
\[L_2+\frac15(\gamma+1)^2+2(\gamma+1)E_2+\frac{\gamma+1}{270}(1440 \log A-108 \gamma-1440 \zeta'(-3)-157)+U+V.\]
\end{thm}

\textit{Proof.} We start with the evaluation of the first integral, and we use formula (4.5) of \cite{CSA} which states that
\[\log\frac{G(1+x)}{G(1-x)}=-x\log\frac{\sin(\pi x)}{\pi}-\frac{1}{2\pi}\Cl_2(2\pi x).\]
By the basic recursion \eqref{GzG} we can easily get the following expression
\[\log^2\frac{G(x)}{G(1-x)}=\lGs(x)+x^2\log^2\frac{\sin(\pi x)}{\pi}+\frac{1}{4\pi^2}\Cl_2^2(2\pi x)+2x\lG(x)\log\frac{\sin(\pi x)}{\pi}+\]
\[\frac{1}{\pi}\lG(x)\Cl_2(2\pi x)+\frac{x}{\pi}\log\frac{\sin(\pi x)}{\pi}\Cl_2(2\pi x).\]
Integrating on $x\in[0,1]$ we have that
\[\int_0^1\log^2\frac{G(x)}{G(1-x)}dx=I_1+\cdots+I_6.\]
We evaluate these six integrals one by one. The first we already know:
\[I_1=L_2.\]
The second integral is
\[I_2=\int_0^1x^2\log^2\frac{\sin(\pi x)}{\pi}dx=\]
\[\frac13\log^2\pi+\int_0^1x^2\log^2(\sin(\pi x))dx-2\log\pi\int_0^1x^2\log(\sin(\pi x))dx.\]
The appearing two integrals are not elementary either; the first -- after an integral transformation -- can be evaluated by a result of Borwein and Straub \cite[Example 5.]{BS}:
\[\int_0^1x^2\log^2(\sin(\pi x))dx=\frac{13\pi^2}{360}+\frac13\log^2(2)+\frac{1}{\pi^2}\log(2)\zeta(3),\]
while the second one can be found in \cite[(5.7)]{EM2}:
\[\int_0^1x^2\log(\sin(\pi x))dx=-\frac13\log2-\frac{\zeta(3)}{2\pi^2}.\]
Hence
\[I_2=\frac{13\pi^2}{360}+\frac13\log^2(2\pi)+\frac{1}{\pi^2}\log(2\pi)\zeta(3).\]
The value of $I_3$ can be given by Parseval's identity:
\[I_3=\frac{1}{4\pi^2}\int_0^1\Cl_2^2(2\pi x)dx=\frac{\pi^2}{720}.\]
In the integral $I_4$ appears the sum of special interest we denoted by $\zeta_H'(2)$.
\[I_4=2\int_0^1x\lG(x)\log\frac{\sin(\pi x)}{\pi}dx=\]
\[2\int_0^1x\lG(x)\log\sin(\pi x)dx-2\log\pi\int_0^1x\lG(x)dx=\]
\[-2\log\pi\left(\frac{\zeta'(2)}{2\pi^2}+\frac16\log(2\pi)-\frac{\gamma}{12}\right)+2\int_0^1x\lG(x)\log\sin(\pi x)dx.\]
Here the exact value of $E_1=\int_0^1x\lG(x)dx$ was used \cite{EM1}. In the literature \cite[Example 7.3]{EM1} we could find
\[\int_0^1\lG(x)\log\sin(\pi x)dx=-\frac12\log2\log(2\pi)-\frac{\pi^2}{24},\]
but not
\[\int_0^1x\lG(x)\log\sin(\pi x)dx,\]
what we would need. Fortunately the generalized Parseval identity \eqref{Parsgen} combined with Theorem \ref{xlogGammaFourier} helps. Recalling that for $\log\sin(\pi x)$ the Fourier sine-coefficients are all zero, and the cosine-coefficients are
\[a_0=-\log2,\quad a_n=-\frac{1}{2n}\;\;(n>0),\]
we have that
\[\int_0^1x\lG(x)\log\sin(\pi x)dx=\]
%\[-\frac14\log2\log\left(\frac{2\pi}{A^4}\right)-\sum_{n=1}^\infty\frac{1}{2n}\left(\frac{1}{4n}-\frac{C}{\pi^2n^2}-\frac{\log n}{4\pi^2n^2}-\frac{1}{\pi^2}T_n\right)=\]
\[-\frac14\log2\log\left(\frac{2\pi}{A^4}\right)-\frac{\pi^2}{48}+\frac{C\zeta(3)}{2\pi^2}+\frac{1}{2\pi^2}\zeta'(3)+\frac{1}{2\pi^2}\sum_{n=1}^\infty\frac{H_n\log n}{n^2},\]
where we used Lemma \ref{Lem1}. Hence
\[I_4=\frac{1}{\pi^2}\zeta_H'(2)-\frac{\pi^2}{24}-\frac12\log(2\pi)\log\left(\frac{2\pi}{A^4}\right)+\frac{\gamma\zeta(3)}{\pi^2}+\frac{1}{\pi^2}\zeta(3)\log(2\pi)+\frac{\zeta'(3)}{\pi^2}.\]
In the determination of $I_5$ there is no additional difficulty:
\[I_5=\frac{1}{\pi}\int_0^1\lG(x)\Cl_2(2\pi x)dx=\frac{1}{2\pi^2}(C\zeta(3)-\zeta'(3)).\]

Lemma \ref{LemCl2} helps us to find
\[I_6=\int_0^1\frac{x}{\pi}\log\frac{\sin(\pi x)}{\pi}\Cl_2(2\pi x)dx=\frac{\zeta(3)\log(2\pi)}{2\pi^2}+\frac{\pi^2}{720}.\]
The expression under \eqref{GperG} is now can be completed by adding the six integrals and making the possible simplifications.

The second statement of the theorem can be established via the expression \cite[eq. (5.5)]{CSA}
\[\log(G(x+1)G(1-x))=-(\gamma+1)x^2-Z(x).\]
Then using that $\log(G(x+1)G(1-x))=\lG(x)+\log(G(x)G(1-x))$, squaring and calculating the resulting integrals one by one, it is possible to prove the second statement of the theorem. One sum arises in the integral of $x^2Z(x)$, and this is
\[\int_0^1x^2Z(x)dx=\sum_{n=1}^\infty\frac{\zeta(2n+1)}{(n+1)(2n+5)}=\frac{1}{540}(1440\log A-1440\zeta'(3)-108\gamma-157),\]
as one can deduce from \cite[p. 515, (2.14)]{CCS}.
\hfill\qed

Now employing \eqref{G2expr} we arrive at our main achievement, the wanted expression for $G_2$.

\begin{thm}We have that
\[4G_2=\log A\left(\frac{10 }{3}+\frac{22}{3} \gamma +6\log(2\pi)\right)+\frac{\zeta ''(2)}{\pi ^2}-\frac{16}{3}\zeta'(-3)(1+\gamma )+\]
\[\frac{\zeta'(3)}{2 \pi ^2}+\frac{1}{\pi^2}\zeta(3)\left(2 \gamma +\frac{1}{2}+3 \log\pi +3\log2\right)+\]
\[\frac{7 \pi ^2}{180}-\frac{11 \gamma ^2}{30}-\frac{157 \gamma }{270}-\frac{103}{270}+\frac{1}{6} \log ^2(2 \pi )+\frac{1}{3} \log(2 \pi )+U+V+\frac{\zeta'_H(2)}{\pi ^2}.\]
\end{thm}

That there are some ``unidentified'' constants in the expression of $G_2$ shows the probable complexity of this quantity.

\section{Some additional remarks}

\subsection{The harmonic zeta function}

In the evaluation of $G_2$ the sum
\[\zeta_H'(2)=\sum_{n=2}^\infty\frac{H_n\log n}{n^2}\]
appeared. The harmonic zeta function, studied first by Apostol and Vu \cite{AV}, is defined as
\[\zeta_H(s)=\sum_{n=1}^\infty\frac{H_n}{n^s}\quad(\Re(s)>1).\]
(Actually, they defined a more general function having the generalized harmonic numbers in the nominator). Hence the notation $\zeta_H'(2)$ already can be understood (apart, of course, from the minus sign left out for the sake of simplicity). Apostol and Vu considered the analytic extension of the $\zeta_H$ function to the whole complex plane, minus discrete singularities described entirely in the paper. No functional equation for this kind of Dirichlet function is known, however, this would be necessary to find a closed form of $\zeta_H'(2)$ (in terms of zeta function values).

This sum appears probably very few times in the mathematical literature. Nevertheless, we could find an expression for it which involves the log-Gamma function:
\[\sum_{n=2}^\infty\frac{H_n\log n}{n^2}=\frac{\pi^2}{6}\left(\frac{\pi^2}{8}-\frac{\gamma^2}{2}+\log^2(2\pi)\right)+\pi^2\gamma\log\left(\frac{A^2}{\sqrt[6]{2\pi}}\right)-C\zeta(3)-\]
\[\zeta'(3)+\frac12\zeta''(2)-2\pi^2\int_0^1x\log^2\Gamma(x)dx,\]
where, as everywhere in the paper, $C=\log(2\pi)+\gamma$, and $A$ is the Glaisher-Kinkelin constant. The proof can be established by using the generalized Parseval identity with $f(x)=x\lG(x)$ and $g(x)=\lG(x)$.

\subsection{An integral of $e^x\lG(x)$}

Before going forward to the evaluation of $G_2$ we remark some corollaries of our results.

Considering the Fourier coefficients $a_1$ and $b_1$ of $x\lG(x)$ we can see that
\[\int_0^1x\lG(x)\sin(2\pi x)dx=\frac{\gamma}{4\pi}.\]
\[\int_0^1x\lG(x)\cos(2\pi x)dx=\frac18-\frac{C}{2\pi^2}-\frac{1}{2\pi^2}\sum_{n=1}^\infty\frac{\log n}{n^2-1}.\]

During our work we have found another integral evaluation:
\[\frac{2}{e-1}\int_0^1 e^x\lG(x)dx=C\left(3-\coth\left(\frac12\right)\right)+\Re\left(\psi\left(\frac{i}{2\pi}\right)\right)-4\sum_{n=1}^\infty\frac{\log n}{1+4n^2\pi^2},\]
where $\psi=\frac{\Gamma'}{\Gamma}$ is the Digamma function \cite{GR,SCbook}. The proof of this integral can be established by the Fourier expansion of $e^x$ together with \eqref{Parsgen}. We note that
\begin{equation}
E_n=\int_0^1x^n\lG(x)dx\label{Endef}
\end{equation}
can be evaluated for any $n\ge0$ via zeta values \cite{EM1}. However, we cannot get our above integral by dividing $E_n$ with $n!$ and summing over $n$ because of the complex structure of the expression for $E_n$.

\subsection{The asymptotics of $G_n$}

We mentioned in the first section that Amdeberhan and his coworkers \cite{ACEKMM} made the observation
\[\frac{L_n}{n!}\sim1.\]
Their approach (considering the Taylor expansion of $\lG(x)$ around $x=0$) can be used to present the asymptotics of $G_n$.

\begin{thm}We have that asymptotically
\[G_n\sim(-1)^nn!.\]
\end{thm}

It is worth to note that in \cite{ACEKMM} a finer description is given for $L_n/n!$, namely,
\[\frac{L_n}{n!}\sim\sum_{i=1}^\infty\frac{(-1)^{i+1}c_i}{i^n}\]
as $n\to\infty$. Here the coefficients $c_i$ can be determined. It would be interesting to find a similar asymptotic expansion for $(-1)^nG_n/n!$.

%\subsection{On the complexity of infinite sums}

%Finally, we would like to give some remarks on the the infinite sum $\zeta_H'(2)$ appeared in the paper.

%It is a well known fact that some integrals can be elementarily evaluated but a tiny modification in the integrand can result very hardly evaluable or even ``inevaluable''\footnote{Inevaluable here means that the integral in question cannot be traced back to more elementary expressions, as it is for some elliptic integrals, for example. See \cite{BM} for a huge collection of beautiful integrals} integrals. The same is true for infinite sums. The complexity of a sum is very hard (or maybe impossible) to measure. A possible insight into the "deeper structure" of an infinite sum (or, better to say, the number what the sum do represent) might be an alternative expression we can find for it. To provide some examples, we know that
%\[e=\sum_{n=0}^\infty\frac1{n!}.\]
%There are a plenty of alternative expressions and appearances for $e$ and it is probably not considered to be a ``hard number''. However, about the modified sum
%\[\sum_{n=1}^\infty\frac1{n\cdot n!},\]
%which appears in the theory of the exponential integral, is harder to see many appearances. Hence, we might say that the complexity of a number is correlated to the quantity of information we know about it.

\section*{Acknowledgement}

The author is grateful to Mr. Donal Connon for his valuable comments and discussions on my work.

\end{document}